\documentclass[reqno,12pt]{article}
\usepackage{amsmath, latexsym, amsfonts, amssymb, amsthm, amscd,graphicx,color,epsfig}
\usepackage{graphics,epsf,psfrag}

\usepackage{amsmath, latexsym, amsfonts, amssymb, amsthm, amscd,graphicx,color,epsfig}
 \usepackage[english]{babel}
\usepackage{graphics,epsf,psfrag}

 \usepackage[utf8x]{inputenc}

\numberwithin{equation}{section}

  \newtheorem{theorem}{Theorem}[section]
  \newtheorem{lemma}[theorem]{Lemma}
  \newtheorem{proposition}[theorem]{Proposition}

\newcommand{\ind}{\mathbf{1}}

\newcommand{\R}{\mathbb{R}}
\newcommand{\Z}{\mathbb{Z}}
\newcommand{\N}{\mathbb{N}}
\renewcommand{\tilde}{\widetilde}
\renewcommand{\hat}{\widehat}

\newcommand{\cC}{{\ensuremath{\mathcal C}} }

\newcommand{\bP}{{\ensuremath{\mathbf P}} }
\newcommand{\bE}{{\ensuremath{\mathbf E}} }

\newcommand{\bN}{{\ensuremath{\mathbf N}} }

\DeclareMathSymbol{\leqslant}{\mathalpha}{AMSa}{"36} 
\DeclareMathSymbol{\geqslant}{\mathalpha}{AMSa}{"3E} 
\DeclareMathSymbol{\eset}{\mathalpha}{AMSb}{"3F}     
\renewcommand{\leq}{\;\leqslant\;}                   
\renewcommand{\geq}{\;\geqslant\;}                   

\newcommand{\ga}{\alpha}
\newcommand{\gb}{\beta}
\newcommand{\gd}{\delta}
\newcommand{\gep}{\varepsilon}       

\newcommand{\gD}{\Delta}
\newcommand{\gk}{\kappa}
\newcommand{\go}{\omega}

\newcommand{\gl}{\lambda}

\newcommand{\gs}{\sigma}

\makeatletter
\def\captionfont@{\footnotesize}
\def\captionheadfont@{\scshape}

\long\def\@makecaption#1#2{%
  \vspace{2mm}
  \setbox\@tempboxa\vbox{\color@setgroup
    \advance\hsize-6pc\noindent
    \captionfont@\captionheadfont@#1\@xp\@ifnotempty\@xp
        {\@cdr#2\@nil}{.\captionfont@\upshape\enspace#2}%
    \unskip\kern-6pc\par
    \global\setbox\@ne\lastbox\color@endgroup}%
  \ifhbox\@ne 
    \setbox\@ne\hbox{\unhbox\@ne\unskip\unskip\unpenalty\unkern}%
  \fi
  \ifdim\wd\@tempboxa=\z@ 
    \setbox\@ne\hbox to\columnwidth{\hss\kern-6pc\box\@ne\hss}%
  \else 
    \setbox\@ne\vbox{\unvbox\@tempboxa\parskip\z@skip
        \noindent\unhbox\@ne\advance\hsize-6pc\par}%
\fi
  \ifnum\@tempcnta<64 
    \addvspace\abovecaptionskip
    \moveright 3pc\box\@ne
  \else 
    \moveright 3pc\box\@ne
    \nobreak
    \vskip\belowcaptionskip
  \fi
\relax
}
\makeatother
\def\writefig#1 #2 #3 {\rlap{\kern #1 truecm
\raise #2 truecm \hbox{#3}}}

\def\beq{\begin{equation}}
\def\eeq{\end{equation}}

\newcommand{\1}{{\rm 1}\kern-0.24em{\rm I}}

\def\_{\vspace{3mm}}

 \title{ A functional limit convergence towards brownian excursion.}
  \author{Julien Sohier \thanks{Universit\'{e} Paris-Dauphine, Ceremade, CNRS UMR 7534, 
   F 75016 Paris France.  e-mail: jusohier@gmail.com } }

\begin{document} 
 
 \maketitle 
 
 \begin{abstract}
   We consider a random walk $S$ in the domain of attraction of a standard normal law $Z$, \textit{ie} there exists a positive 
 sequence $a_n$ such that $S_n/a_n$ converges in law towards $Z$. The main result of this note is that the rescaled process 
 $(S_{\lfloor nt \rfloor}/a_n, t \geq 0)$
 conditioned to stay non-negative, to start and to come back \textit{near the origin}
  converges in law towards the normalized brownian excursion.  
 \end{abstract}

  
 
 \section{ Introduction and the main result }

  It is a classical result that if a random walk $S$ is in the domain of attraction of the standard normal law with 
    norming sequence $a_n$, the rescaled process $(S_{\lfloor nt \rfloor/a_n})_{t \geq 0}$ converges in law towards the 
 brownian motion (see \cite{Bil}). Denoting by $S^{\ast},\bP_x$ the random walk starting from $x$ and conditioned to 
 stay always positive (one can make sense of such a definition by means of a so called $h$-transform), it has
 recently been shown in \cite{BJD} and in \cite{CCa} that if $ x/a_n$ vanishes as $n \to \infty$, 
 the corresponding rescaled process converges in law towards the brownian meander. A natural question related to these results is  
 whether conditioning on a \textit{late} return near the origin (ie on $\{ S^{\ast}_n = y \}$ with $y/a_n \to 0$ as $n \to \infty$)
 implies the convergence of $(S^{\ast},\bP_x)$ towards the brownian excursion. 

  Extending previous results from \cite{BJD}, we show in this paper that such a convergence holds. 
 Before stating precisely our main results, we recall the essentials of the conditioning to stay positive for
  an oscillating random walk. 

 \subsection{ Conditioning a random walk to stay positive} 
  Let $S_n = X_1 + \ldots + X_n$ be an integer valued aperiodic random walk.  We 
 write $\bP_x$ the law of $S$ started at $x$ and for convenience we put $\bP = \bP_0$. 
   
  Next we introduce the strict descending  \textit{ladder process} $(T^-_k,H^-_k )_{k \geq 0}$
  by setting $ (T^-_0, H^-_0) = (0,0)$ and 
   \begin{equation}
     T^-_{k+1} := \min \{ j >  T^-_k | S_j < S_{T^-_k} \}, \hspace{ 1 cm} H^-_{k+1} = -S_{T^-_{k+1}}.
   \end{equation}
 Note that under $\bP$, $(T^-, H^-)$ is a bivariate renewal process, that is a random walk on $(\Z^+)^2$ with step law supported on 
the first quadrant. The sequence $T^-$ is the sequence of the so called \textit{ (strictly) descending
 ladder epochs}, the sequence $H^-$ the sequence of 
 \textit{descending ladder heights}. 

    We denote by $V(\cdot)$ the renewal function associated to $H^-$, that is the positive function defined by
   \begin{equation}
   V(x) := \sum_{k \geq 0} \bP( H^-_k \leq x). 
 \end{equation}
  Note in particular that $V(y)$ is the expected number of ladder points in the stripe $[0,\infty) \times [0,y]$. It follows
 that it is a subadditive and increasing function. 

  The \textit{killed} random walk $\hat{S}$ is a Markov chain defined in the following way.
 Let $\tau_{(-\infty,0)}$ denote the first entrance time of $S$ into the negative half plane.
 Introducing $\{ \gD \}$ a cimetery state, for every $n$, 
 \begin{equation}
   \hat{S}_n := S_n \ind_{ \tau_{(-\infty,0)} > n} + \gD \ind_{ \tau_{(-\infty,0)} \leq n}.
 \end{equation}
 
    Then we denote $S$ conditioned to stay non negative by $S^{\ast}_n = \sum_{i = 1}^n X^{\ast}_i$. In
 our integer valued oscillating case
 this is a Markov chain on $\Z^+$ whose law is defined for any $n \in \N$ and for any $B \in \gs(S_1, \ldots, S_n)$ by: 
   \begin{equation}
   \bP^{\ast}_x[ B \cap \{ S_n = y \}] := \frac{ V(y)}{V(x)} \bP_x[ B \cap \{ S_n = y \} \cap \mathcal{C}_n ]
 = \frac{V(y)}{V(x)}\bP_x[ B \cap \{ \hat{S}_n = y \} ],
 \end{equation}
 where $\mathcal{C}_n = \{ S_1 \geq 0, \ldots, S_n \geq 0 \}$. The terminology is justified by the following 
 weak convergence result
 \begin{equation}
   \bP^{\ast}_x   = \lim_{n \to \infty} \bP_x \left( \cdot   | \mathcal{C}_n \right)
 \end{equation}
  which is proved in \cite{BD}, Theorem 1. 
  
 \subsection{ A convergence towards the brownian excursion}
   From now on,  we will always assume that
 $S$ lies in the domain of attraction of the standard normal law.
  This means that the sequence $(X_k)$ is iid and that for a suitable norming sequence $(a_n)$ one has the weak convergence
  \begin{equation} \label{cd}
    S_n/a_n \Rightarrow \phi(x) dx , \hspace{.6 cm} \phi(x) := \frac{1}{\sqrt{2 \pi}} e^{-x^2/2}.
  \end{equation}
  In particular this is the case when $ \bE[X_1] = 0$ and $\bE[X_1^2 ] =: \gs^2 < \infty$ with $a_n = \gs \sqrt{n}$ by the central limit 
 theorem.
  
   By standard theory of stability, (see \cite{Fel2} IX.8 and XVII.5) for \eqref{cd} to hold it is necessary and sufficient that
 $\bE[X_1] = 0$, that the truncated variance $ \Phi(t) := \bE[ X_1^2 \ind_{ |X_1| \leq t}]$ is slowly varying at infinity (that is 
 $ \frac{\Phi(ct)}{\Phi(t)} \to 1 $ as $t \to \infty$ for any $c > 0$ ) and that the sequence $a_n$ satisfies $ a^2_n \sim n \Phi(a_n)$
 as $n \to \infty$. 

  We define $ \Omega$ as being the space $ D([0,1],\R)$ the set of c\`{a}dl\`{a}g functions on $[0,1]$ endowed with the standard Skorohod
 topology (see \cite{Bil}) and for $n \in \Z^+$, we define the application $X^n$ by:
 \begin{equation} 
  X^n: \begin{array}{ccc}
     \Z^n & \longrightarrow &  \Omega\\
    (u_1, \ldots, u_n) & \mapsto & \left( \frac{\sum_{i=1}^{[nt]} u_i }{a_n} \right) _{t \in [0,1]}  
  \end{array}.
   \end{equation}
 For $x,y$ positive integers, we denote by $P_n^{\ast,x,y}$ the law of $S^{\ast}$ 
conditionally on the event $\{  S^{\ast}_0 = x, S^{\ast}_n = y \}$, 
  and we define the probability laws on $\Omega$:
   \begin{equation}
     Q_n^{x,y} := P_n^{\ast,x,y} \circ (X^n)^{-1}.
   \end{equation}
   We can now state our main result:
    \begin{theorem} \label{TA}
     Let $x_n$ and $y_n$ be positive integer valued sequences such that $x_n/a_n \to 0$ and $y_n/a_n \to 0$.  
 Then, as $ n \to \infty$, the following convergence holds in $\Omega$:
 \begin{equation}
   Q_n^{x_n,y_n} \Rightarrow e
 \end{equation}
   where $e$ denotes the law of the normalized brownian excursion. 
\end{theorem}
  
 The proof of this result will follow the standard procedure of showing finite dimensional convergence and tightness.

   \subsection{ Some motivations and a short overview of the literature } The study of invariance principles for 
 random walks is a very classical topic in probability (classical references are \cite{Sko}, \cite{Bil}).
  Extending these invariance principles
 to conditioned random walks is far from being straightforward. Sometimes a clever representation can considerably 
 simplify the proofs (like in \cite{Bolt}, \cite{Do2} for the convergence towards the meander),
  but generally speaking such an issue demands some technical efforts, see \cite{Igl} for a convergence towards the meander or
 \cite{Lig} for the brownian bridge. 

   The more particular case of convergence towards the brownian excursion for the conditioned simple random walk conditioned 
 by a late return to zero has first been proved in \cite{DIM}. Their results have been extended to the case where $S$ has finite variance
 in \cite{K}. 

    A related result to ours that will turn out to be quite useful in our proofs is  the convergence towards the brownian meander 
 of a random walk in the domain of attraction of the normal law  starting from $x_n$ where $x_n$ is $o(a_n)$ conditioned on $\mathcal{C}_n$ 
 (see  \cite[Remark 4]{Shi}). Combining tightness arguments and local limit estimates, 
 this result has been extended to the case where $S$ is conditioned to stay positive by \cite{BJD}, and 
 their results in turn have been extended by quite different and somewhat lighter techniques in \cite{CCa}
 to the case where $S$ is in the 
  domain of attraction of a \textit{stable law} 
 with index $\ga \in (0,2]$ and with positivity parameter $\rho \in (0,1)$. Lacking a
 suitable representation under the form of an $h$-transform
  for the brownian excursion, our methods follow the same path as in \cite{BJD}. 

  Besides the interest they have in their own, invariance principles are important in view of their applications.  Let us 
 mention one of them which is actually the main motivation of this paper. Consider the following homogeneous polymer model (a by now
 classical reference for polymer models is \cite{GB}): for $N \in \N, y \in \R^+, a > 0
$ and $\gep \in \R$, we set
 \begin{equation}  
    \frac{d \bP_{N,a,\gep}^c}{d \bP} := \frac{1}{ Z_{N,a,\gep}} \exp \left( \gep \sum_{i=1}^N  \ind_{ S_i \in [0,a]} \right) \ind_{ S_N \in [0,a]}        
       \end{equation}
 where $\bP$ is an aperiodic $\Z$ valued random walk in the domain of attraction of the standard normal law. The law $\bP_{N,a,\gep}^c$ 
 may be viewed as an effective model for a $(1+1)$ dimensional interface above a wall with homogeneous impurities which are concentrated
 in the stripe $[0,\infty) \times [0,a]$. These impurities are either attracting or repelling the interface (depending on the sign of $\gep$). 
 
 One standard goal related to this kind of models is to find the asymptotic behavior of the typical paths in the limit $N \to \infty$ 
 and to study their dependence on $\gep$ and $a$. These limits have been resolved in the thesis \cite{JS}.

  A common feature shared by this model and the classical homogeneous one is that the measure $ \bP_{N,a,\gep}^c$ exhibits
 a remarkable decoupling between the contact level set $\mathcal{I}_N := \{ i \leq N, S_i \in [0,a] \}$ and the excursions of $S$ between 
 two consecutive contact points (see \cite{DGZ} for more details in the standard homogeneous pinning case). In fact, conditionally on 
  $I_N = \{ t_1, \ldots, t_k\}$ and on $(S_{t_1}, \ldots, S_{t_k})$, the \textit{bulk} excursions 
 $e_i = \{ e_i(n) \}_n := \left\{ \{S_{t_i+n} \}_{ 0 \leq n \leq t_{i+1} - t_i} \right\}$
 are independent under  $ \bP_{N,a,\gep}^c$ and are distributed like the random walk $(S,\bP_{ S_{t_i}})$ conditioned on the 
 event  $\left\{ S_{t_{i+1} - t_i} \in [0,a], S_{t_i + j} > a, j \in \{ 1, \ldots, t_{i+1} - t_i -1 \} \right\}$. It is therefore clear that 
 to extract scaling limits on $ \bP_{N,a,\gep}^c$, one has to combine good control over the law of the contact set $\mathcal{I}_N$ 
 and suitable asymptotics properties of the excursions, and for this the utility of Theorem \ref{TA} emerges (see chapter 3 of
 the thesis \cite{JS} for details).

 \subsection{Outline of the paper} The exposition of this paper will be organized as follows:
 \begin{enumerate}
  \item[-]  in Section \ref{SecPre}, we collect some preliminary facts.
  \item[-]  in Section \ref{SecFDC}, we discuss finite dimensional convergence and state our main technical lemma.
  \item[-] in Section \ref{MMM}, we prove Lemma \ref{KL}, which implies the finite dimensional convergence in Theorem \ref{TA}.
  \item[-] in Section \ref{SecTI}, we show the tightness of the sequence of measures $ ( Q_n^{x_n,y_n})_n$, thus proving Theorem \ref{TA}.
  \item[-] in Section \ref{SecV}, we give a uniform equivalence for the tails of the random variable $\tau_{(-\infty,0)}$ under the
  law $\bP_{x_n}$. This estimate 
 is widely used in sections \ref{MMM} and \ref{SecTI}. 
 \end{enumerate}

 \section{Some preliminary facts}  \label{SecPre}
 
  \subsection{ Regular varying sequences}  Throughout this note, for positive sequences $\ga_n$ and $\gb_n$,
 we use the notation $\ga_n \sim \gb_n$ 
 to indicate that $ \ga_n/\gb_n \to 1$ as $n \to \infty$. Following Doney's terminology,
 for positive measurable functions $g,h$ on $\R^+$, 
  we will often say that the equivalence 
  \begin{equation}
    g(x_n) \sim h(x_n)
  \end{equation} 
 is true 
 \textit{uniformly on the sequences $x_n$ such that $x_n/a_n \to 0$}. By this we mean that, given
 any positive sequence $\gep_n$ such that
 $\gep_n \to 0$ as $n \to \infty$, the convergence 
  \begin{equation}
    \frac{ g(x_n)}{h(x_n)} \to 1
  \end{equation}
 holds uniformly for every sequence $x_n \in \Delta_{ \gep_n}$ where 
 \begin{equation}
  \Delta_{\gep_n} := \{ y \in \Z^{\N}, \forall n \geq 0, y_n \in [0, \gep_n a_n] \}.
 \end{equation}

 A positive sequence $d_n$ is said to be slowly varying with index $\ga \in \R$ (which we denote by 
 $d_n \in \R_{\ga}$) if $d_n \sim L_n n^{\ga}$ where $L_n$ is slowly varying at infinity that is for every positive $t$, 
 $\lim_{n \to \infty} \frac{L_{[nt]}}{L_n} = 1$. If $d_n \in \R_{\ga}$, we can (and will always assume) that $d_n = d(n)$ where 
 $d(\cdot)$ is a continuous strictly monotone function whose inverse will be denoted $d^{-1}(\cdot)$ (see \cite[Theorem 1.5.3]{BGT}). 
  Observe that if $d_n \in \mathcal{R}_{\ga}, 
 d^{-1}(n) \in \mathcal{R}_{1/\ga}$ and $1/d_n \in \mathcal{R}_{- \ga}$. 

   The following basic uniform convergence property (\cite[Theorem 1.2.1]{BGT}) will be often used in the sequel;
  if $d_n \in \mathcal{R}_{\ga}$, then for every fixed $\gep > 0$ 
 \begin{equation} \label{RVV}
   d_{[tn]} = t^{\ga}d_n (1 + o(1))  
 \end{equation}
 uniformly for $t \in [ \gep, 1/\gep]$. \\

 \subsection{ Fluctuation theory} \label{FlT}
 In a similar way as for the descending ladder process,
  one can define the weak ascending bivariate renewal process $(T_k^+,H_k^+)_k$ as $T^+_0 := 0$,
 $T_{k+1} := \min \{ j > T_k^+, S_j \geq S_{T_k^+} \}$, 
 $H_k^+ := S_{T_k^+}$ and 
 \begin{equation}
         U(x) := \sum_{k \geq 0} \bP( H^+_k \leq x).
      \end{equation}
  It is known that $S_1$ is in the domain of attraction (without centering) of a stable law if and only if $(T^-_1, H_1^-)$ lies 
 in a bivariate domain of attraction (see for example \cite{DGr}). We can specialize this fact to our setting.
  By hypothesis, $S_1$ lies in the domain of attraction of 
 the standard normal law, so that by standard fluctuation theory, $a_n \in \mathcal{R}_{1/2}$. We then define two sequences
 \begin{equation}
   \log( \frac{n}{\sqrt{2}}) = \sum_{m = 1 }^{\infty} \frac{\bP[ S_m < 0]}{m} e^{ - \frac{m}{b_n}}, \hspace{.6 cm} c_n := a(b_n).
 \end{equation}
  Then $b_n \in \mathcal{R}_2$, $c_n \in \mathcal{R}_1$ and we have the weak convergence
 \begin{equation} \label{BI}
   \left( \frac{T^-_n}{b_n}, \frac{H^-_n}{a_n} \right) \Rightarrow Z, \hspace{.6 cm} 
 \bP[ Z \in (dx,dy) ] = \frac{e^{-1/2x}}{\sqrt{2 \pi} x^{3/2}} \ind_{ x \geq 0} \gd_1(dy),
 \end{equation}
  where $\gd_1(dy)$ denotes the Dirac measure at $y = 1$.  Note in particular that, like in the simple random walk case,  $T_1^-$ is 
 attracted to $Y$, the stable law of index $1/2$.
  \begin{equation}
    \frac{T_n^-}{b_n} \Rightarrow Y, \hspace{.6 cm} \bP[ Y \in dx] = \frac{e^{-1/2x}}{\sqrt{2 \pi} x^{3/2}} \ind_{ x \geq 0}.
  \end{equation}
  We recall also that $b_n$ is sharply linked to the tails of $T_1^-$ by the relation
 \begin{equation}
   \bP[ T_1^- > b_n ] \sim \sqrt{\frac{2}{\pi}} \frac{1}{n} 
 \end{equation}
 and it is known that this is a necessary and sufficient relation in order for $b_n$ to be such that $ T_n^-/b_n \Rightarrow Y$. 
    
   Equation \eqref{BI} also implies that the process $(H^-)$ follows a 
 generalized law of large numbers, namely
 $ \frac{H^-_n}{c_n} \Rightarrow 1$ ($H^-_1$ is said to be \textit{relatively stable}).
  Consequently the following 
 equivalence holds (see \cite[Theorem 8.8.1]{BGT}) 
 \begin{equation} \label{UVeq}
   V(x) \sim c^{-1}(x) =: \frac{ x}{l^-(x)} 
 \end{equation}
 where $l^-(\cdot) $ is  slowly varying at infinity. In a similar way, one can prove that the equivalence
  \begin{equation}
     U(x) \sim \frac{ x}{l^+(x)} 
  \end{equation}
  is verified for some slowly varying function $l^+(\cdot)$.

  \subsection{ The duality lemma and local limit estimates} 
   Let  $v(\cdot,\cdot)$ be the renewal
  mass function of the bivariate renewal process $(H^-,T^-)$, that is 
 \begin{equation}
  v(n,x) := \sum_k \bP[ T_k^- = n, H_k^- = x]
 \end{equation}
   and $u(\cdot,\cdot)$ its counterpart for the process $(H^+,T^+)$
   \begin{equation}
  u(n,x) := \sum_k \bP[ T_k^+ = n, H_k^+ = x].
    \end{equation}
   The power of fluctuation theory for the study of random walks is linked to some fundamental identities, the most famous one 
 being the so called "duality lemma" (see \cite[Chapter XII]{Fel2} ):
 \begin{equation}
   \bP[ \mathcal{C}_n, S_n \in dx] = \bP [ n \hspace{.05 cm} \text{ is a ladder epoch}, S_n \in dx] = u(n,x)
 \end{equation}
 where by the event $\{  n \hspace{.05 cm} \text{ is a ladder epoch} \}$ we mean of course the disjoint union of the events 
 $\cup_k \{ T_k^+ = n \}$. 
  The following equivalence about the asymptotics of $u(\cdot,\cdot)$ has been shown independently in \cite{Car} and in \cite{BJD}.
 Note that for the later, it is the chore of 
 the proof of their main result. 
  \begin{lemma} \label{BJD1}
    Uniformly for $ 0 \leq y_n \leq K a_n$, one has the following equivalence:
  \begin{equation}
    \bP[\hat{S}_n = y_n] = u(n,y_n) \sim \frac{U(y_n)}{n} \bP[S_n = y_n].
  \end{equation}
 \end{lemma} 

  \section{ Finite dimensional convergence in Theorem \ref{TA}} \label{SecFDC} 
 
   \subsection{ The law of the renormalized brownian excursion}  
   For $x,y,t > 0$, we define $q_t(x,y)$ the transition function of the killed Brownian motion, that is 
 \begin{equation}
   \begin{split}
      &  q_t(x,y) := \frac{1}{\sqrt{t}} r( \frac{x}{\sqrt{t}},\frac{y}{\sqrt{t}}) \\
   \text{where} \hspace{.2 cm} & r(u,v) := \sqrt{\frac{2}{\pi}}  \sinh(uv) \exp(-\frac{u^2 + v^2}{2}), 
   \end{split}
 \end{equation}
   and the following transition function :   
 \begin{equation}
    \begin{split}
     & l_t(y) := \frac{1}{t} r_0( \frac{y}{\sqrt{t}}) \\
     \text{where} \hspace{.2 cm}  & r_0(v) := \sqrt{\frac{1}{2\pi}} v  \exp(-\frac{v^2}{2}).
    \end{split}
  \end{equation}
    It is well known that (see \cite{Bo}) for $ k \in \bN, 0 < t_1 < \ldots < t_k < 1$ and $f \in \cC^b([0,1]^k,\R)$, one has:
 \begin{equation}
\begin{split}
  & e(f(\go_{t_1}, \ldots, \go_{t_1})) \\
  & \phantom{iiiiiiiii} = 2 \sqrt{ 2\pi} \int_{(\R^+)^k} f(x_1,\ldots,x_k) l_{t_1}(x_1) \ldots
 q_{t_k-t_{k-1}}(x_{k-1},x_k)l_{1-t_k}(x_k) dx_1\ldots dx_k  . 
\end{split}
 \end{equation}

 To get Theorem \ref{TA}, we have to show finite dimensional convergence, that is we show that for every positive integer $k$, 
  $(t_1, \ldots, t_k) \in (0,1)^k, f \in \cC^b((\R^+)^k,\R)$:
  \begin{equation} \label{AM5}
     \begin{split}
      & \frac{\bE_{x_n}[ f(\frac{S^{\ast}_{\lceil nt_1 \rceil}}{a_n}, \ldots, \frac{S^{\ast}_{\lceil nt_k \rceil}}{a_n}) \ind_{S^{\ast}_n = y_n}]} { \bP_{x_n} [S^{\ast}_n = y_n] } \\
     & \phantom{iiiiii} \to 2 \sqrt{2\pi} \int_{\R^+} f(x_1, \ldots, x_k) l_{t_1}(x_1) q_{t_2-t_1}(x_1,x_2) \ldots 
 l_{1-t_k}(x_k) dx_1 \ldots dx_k 
       \end{split}
   \end{equation}
  as $n \to \infty$.
 
  \subsection{ Getting the convergence \eqref{AM5} }
 Our main tool to get this convergence is the following result which we prove in part \ref{MMM}:
 \begin{lemma} \label{KL}
   For $K > 0$, uniformly in $x_n/a_n \to 0$ as $n \to \infty$ and in $y_n$ such that $y_n/a_n \in [0,K]$,
 one has the following equivalence:
   \begin{equation}
     \bP_{x_n}(\hat{S}_n = y_n) \sim \frac{V(x_n) U(y_n)}{n} \bP(S_n = y_n).
   \end{equation}
   \end{lemma}
       The next result is a consequence of the Wiener Hopf factorization, it has been shown in \cite{BJD} and it will turn out to be
   useful numerous times in the sequel.
         \begin{lemma}\label{BJDo} Let $K > 0$. Uniformly in the sequences $(x_n)_{n \geq 0},  (y_n)_{n \geq 0}$  
    such that $x_n/a_n \in [0,K], y_n/a_n \in [0,K]$, one has the following equivalence: 
      \begin{equation}
       \frac{ U(x_n) V(y_n)}{n} =  2 \frac{x_n}{a_n} \frac{y_n}{a_n} + o(1) \hspace{ .3 cm} \text{as} \hspace{ .3 cm} n \to \infty
        \end{equation}
         \end{lemma}
         Lemma \ref{KL} straightforwardly implies the equivalence:
          \begin{equation} \label{TRE}
           \bP_{x_n} (S^{\ast}_n = y_n) \sim  \frac{ U(y_n) V(y_n)}{n} \bP[ S_n = y_n].
            \end{equation} 
        Of course, $S^{\ast}$ is not reversible. Nevertheless, using time reversal, combining Lemma \ref{BJDo} and
 the equivalence \eqref{TRE}
        straightforwardly imply the following:
         \begin{lemma} \label{DDD}
          For $K > 0$, uniformly in $x_n/a_n \in [0,K]$  and in $y_n$ such that $y_n/a_n \to 0 $ as $n \to \infty$, one 
has the following equivalence:
           \begin{equation}
            \bP_{x_n}(S^{\ast}_n = y_n) \sim 2\frac{y_n^2}{a_n^2} \bP(S_n = x_n) \sim 2\frac{y_n^2}{a_n^2} \frac{ \phi(x_n/a_n)}{a_n}.
            \end{equation}
          \end{lemma}
          We finally recall the following proposition from \cite{BJD}:
          \begin{proposition} \label{BBJ}
             Suppose $x_n$ and $y_n$ are integers such that
          \begin{equation}
           x_n/a_n \to u > 0, \hspace{1 cm} y_n/a_n \to v > 0 
           \end{equation}
            as $n \to \infty$. Then one has the convergence:
              \begin{equation} \label{CUU}
                a_n \bP[\hat{S}_n = y_n ] \to r(u,v)
                \end{equation}   
         \end{proposition}
        
    It is then easy to check that combining the Lemmas \ref{KL}, \ref{BJDo}, \ref{DDD} and the Proposition \ref{BBJ}, 
 one gets the convergence in \eqref{AM5}, so that finite dimensional convergence in Theorem \ref{TA} holds. 
 
  \section{Proof of Lemma \ref{KL}} \label{MMM}

 \subsection{ The case where $y_n/a_n$ is bounded away from zero } 
  
   We first assume that there exists $\gep > 0$ such that for every $n,  y_n/a_n \geq \gep$. 

  We define $m_n := \inf \{S_j,j\leq n \}$ and $\mu_n := \inf \{j \leq n, S_j = m \}$ and their all time counterparts 
  $m = \inf \{S_j,j\geq 0 \}$ and $\mu := \inf \{j \geq 0, S_j = m \}$. Let $\eta > 0$ be fixed.

  Alili and Doney  have used  the following equality in \cite{AD}, it is an easy consequence of the duality lemma:
  \begin{equation} \begin{split}
     \bP_{x_n}[\hat{S}_n = y_n] & =  \bP_{x_n}[\hat{S}_n = y_n ; \mu_n < \eta n] +  \bP_{x_n}[\hat{S}_n = y_n ; \mu_n \geq \eta n] \\
                            & = \sum_{j=0}^{ \eta n} \sum_{k=0}^{x_n \wedge y_n} \bP_{x_n}[S_n = y_n, \mu_n = j, m_n = k] +  \bP_{x_n}[\hat{S}_n = y_n ; \mu_n \geq \eta n] \\
& = \sum_{j=0}^{ \eta n} \sum_{k=0}^{x_n \wedge y_n} v(j,x_n-k) u(n-j,y_n-k) +  \bP_{x_n}[\hat{S}_n = y_n ; \mu_n \geq \eta n].
     \end{split} 
  \end{equation}
   We first treat the first term in the right hand side of the above equality. The assumptions on $x_n, y_n$ imply that for 
large enough $n$, $x_n \wedge y_n = x_n$. Using Lemma \ref{BJD1}, for large enough $n$, we get that:
   \begin{equation} \begin{split}
     \sum_{j=0}^{ \eta n} \sum_{k=0}^{x_n \wedge y_n} v(j,x_n-k) u(n-j,y_n-k)
        & \sim  \sum_{j=0}^{ \eta n} \sum_{k=0}^{x_n}  v(j,x_n-k) \frac{U(y_n-k) \bP_k[ S_{n-j} = y_n]}{n-j} 
       \end{split}
    \end{equation}
   as $n \to \infty$, so that:
   
    \begin{equation} \begin{split}
       g_n(\eta) \sum_{j=0}^{ \eta n} \sum_{k=0}^{x_n}  v(j,k) 
        &  \leq   \frac{n} { U(y_n)  \bP[S_n = y_n]} \sum_{j=0}^{ \eta n} \sum_{k=0}^{x_n}  v(j,x_n-k) \frac{U(y_n-k) \bP_k[ S_{n-j} = y_n]}{n-j} \\
         &  \leq f_n(\eta) \sum_{j=0}^{ \eta n} \sum_{k=0}^{x_n}  v(j,k) 
          \end{split}
     \end{equation} 
      where we defined 
       \begin{equation}
        f_n(\eta) := \sup_{ j \leq \eta n, k \in [0,x_n]} \frac{U(y_n - k) \bP_k[ S_{n-j} = y_n]}{(1- \eta) U(y_n)  \bP[S_n = y_n]}
        \end{equation}
     and 
        \begin{equation}
        g_n(\eta) := \inf_{ j \leq \eta n, k \in [0,x_n]} \frac{U(y_n - k) \bP_k[ S_{n-j} = y_n]}{ U(y_n)  \bP[S_n = y_n]}.
        \end{equation}
     
      Using the standard local limit theorem and equivalence \eqref{UVeq}, one gets easily that
  $\lim_{ \eta \searrow 0}\limsup_{n \to \infty} f_n(\eta) =  \lim_{ \eta \searrow 0}\liminf_{n \to \infty} g_n(\eta) = 1$.
    Thus we are left with showing that
         \begin{equation}\label{fac}
       \sum_{j=0}^{ \eta n} \sum_{k=0}^{x_n}  v(j,k) \sim V(x_n) .
                  \end{equation}
   Note that of course 
    \begin{equation}
       \sum_{j=0}^{ \infty} \sum_{k=0}^{x_n}  v(j,k) = V(x_n) ,
    \end{equation}
   so that we just have to show that 
   \begin{equation} \label{CVU}
     \frac{\sum_{j > \eta n} \sum_{k=0}^{x_n}  v(j,k)}{V(x_n)} \to 0
   \end{equation}
   as $n \to \infty$ uniformly on $x_n$ such that $x_n/a_n \to 0$. For this, we note that Lemma \ref{BJD1} implies 
    \begin{equation} 
      v(n,x) \sim \frac{V(x) \bP[ S_n = -x]}{n} 
    \end{equation}
    as $n \to \infty$ uniformly on $x \in [0,K a_n]$ where $K > 0$ , so that 
   \begin{equation}
     \sum_{j > \eta n} \sum_{k=0}^{x_n}  v(j,k) \sim \sum_{j > \eta n} \sum_{k=0}^{x_n} \frac{ V(k) \bP[ S_j = - k]}{j}.
   \end{equation}
 Using the fact that $V(\cdot)$ is increasing and the standard local limit theorem
  (here and later $c$ is a positive constant which may vary from line to line):
   \begin{equation}
     \frac{\sum_{j > \eta n} \sum_{k=0}^{x_n}  v(j,k)}{V(x_n)} \leq c \sum_{j > \eta n} \sum_{k=0}^{x_n} 
\frac{ \phi( k/a_j)}{ j a_j} \leq c \sum_{j > \eta n} \frac{x_n} { j a_j}.
   \end{equation}
  Finally, as $a_n \in \mathcal{R}_{1/2}$, using property \eqref{RVV} it is easy to see that 
   \begin{equation}
     \sum_{j > \eta n} \frac{a_n} { j a_j} \sim \int_{\eta}^{\infty} x^{-3/2} dx 
   \end{equation}
  and this entails \eqref{fac}.
 To conclude the case where $y_n/a_n$ is bounded away from zero,  we are left with showing that for any $\eta > 0$, one has: 
 \begin{equation}
    \limsup_{ n \to \infty} \frac{ n \bP_{x_n}[\hat{S}_n = y_n ; \mu_n \geq \eta n]}{ V(x_n)U(y_n) \bP[ S_n = y_n]}   = 0
 \end{equation}
 as $n \to \infty$. By the standard local limit theorem, there exists $a,b >0$ such that $a \leq a_n \bP[ S_n = y_n] \leq b $.
  Using Lemma \ref{BJDo},
 we get that: 
 \begin{equation}
   \frac{ n \bP_{x_n}[\hat{S}_n = y_n ; \mu_n \geq \eta n]}{ V(x_n)U(y_n) \bP[ S_n = y_n]} = \frac{ n \bP^{\ast}_{x_n}[S_n = y_n ; \mu_n \geq \eta n]}{ V(y_n)U(y_n) \bP[ S_n = y_n]} \leq \frac{ c a_n  \bP^{\ast}_{x_n}[S_n = y_n ; \mu_n \geq \eta n]}{ \gep^2 }
 \end{equation}
 so that we have to show that 
 \begin{equation} \label{EV}
   \limsup_{ n \to \infty} a_n \bP^{\ast}_{x_n}[S_n = y_n ; \mu_n \geq \eta n] = 0.
 \end{equation}

   Then we fix $\theta \in (\eta,1)$  and we have:
   \begin{equation} \label{2P} 
   \begin{split}
     a_n \bP^{\ast}_{x_n}[S_n = y_n ; \mu_n \geq \eta n] 
     & = \underbrace{a_n \bP^{\ast}_{x_n}[ \eta n \leq  \mu_n \leq \theta n]}_{(1)} \\
     & + \underbrace{a_n \bP^{\ast}_{x_n}[ \mu_n > \theta n]}_{(2)}.
  \end{split}
 \end{equation}

    Making use of the Markov property, one gets:  
   \begin{equation}
   \begin{split}
        (1) & =  a_n \sum_{ j= \eta n}^{ \theta n} \sum_{k=0}^{x_n}  \bP^{\ast}_{x_n}\left[ \mu_n = j, m_n =k\right]
               \bP^{\ast}_k \left[ S_{n-j} = y_n, \min_{ l \leq n-j} S_l \geq k\right] \\ 
             & \leq a_n \sum_{ j= \eta n}^{ \theta n} \sum_{k=0}^{x_n}  \bP^{\ast}_{x_n}[ \mu_n = j, m_n =k]
 \frac{V(y_n)}{V(k)} \bP_k \left[ \hat{S}_{n-j} = y_n, \min_{ i \leq n-j} \hat{S}_i \geq k\right].                                
   \end{split}
   \end{equation}
   Noting that one has the equality $ \bP_k \left[ \hat{S}_{n-j} = y_n, \min_{ i \leq n-j} \hat{S}_i \geq k\right] =
 \bP \left[ \hat{S}_{n-j} = y_n -k \right]$,
 we get (note that $V(k) \geq 1$ for every $k$):
 \begin{equation}
   (1)\leq a_n \sum_{ j= \eta n}^{ \theta n} \sum_{k=0}^{x_n}  \bP^{\ast}_{x_n}\left[ \mu_n = j, m_n =k\right] V(y_n) 
\bP \left[ \hat{S}_{n-j} = y_n -k \right].
 \end{equation}

   Making use of Lemma \ref{BJDo}, of Lemma \ref{BJD1} and of the fact that
  $x_n / a_n \to 0$ as $n \to \infty$, we get : 
  \begin{equation} 
     \begin{split} 
    (1) & \leq c a_n \sum_{ j= \eta n}^{\theta  n} \sum_{k=0}^{x_n} \bP^{\ast}_{x_n} [ \mu_n = j, m_n=k] { U(y_n) V(y_n) \over n-j} \bP [ S_{n-j} = y_n -k ] \\
         & \leq c K^2  \sum_{ j= \eta n}^{ \theta  n} \bP^{\ast}_{x_n}[ \mu_n = j] { n \over n-j} a_n \bP [ S_{n-j} = y_n].
   \end{split} 
  \end{equation}

 Making use of the standard local limit theorem, we have easily:
  \begin{equation} 
   (1)  \leq  c K^2 (1- \theta)^{ -3/2} \bP^{\ast}_{x_n}\left[ \mu_n \geq \eta n\right].
   \end{equation}
   Evidently, for every $n$, one has $ \mu_n \leq \mu$, so that 
   \begin{equation}
     \bP^{\ast}_{x_n}[ \mu_n \geq \eta n] \leq \bP^{\ast}_{x_n}[ \mu \geq \eta n],
   \end{equation}
   and it has been shown in \cite[Theorem 5.1]{BJD} that for every $\eta > 0$, uniformly in the sequences $x_n$ such that
  $x_n/a_n \to 0$ as $n \to \infty$, 
 the quantity $ \bP^{\ast}_{x_n}[ \mu \geq \eta n]$
 vanishes as $n \to \infty$. 

 For the second term in \eqref{2P}, we will need the following result which has been proved in  \cite{BJD}:
 
  \begin{proposition} \label{BBBB}
    For any $\gk > 0$, for $x_n/a_n \to 0$ as $n \to \infty$, one has the following convergence:
   \begin{equation}
    \bP^{\ast}_{x_n} \left[ \max_{ j \leq \mu} S_j \geq \gk a_n \right] \to 0.
    \end{equation}
   \end{proposition}
    
     We give us $ \gk \in (0, \gep)$ and for $n > 0$, we note $\tau := \inf \{ j \geq 0, S_j \geq \gk a_n \}$. Then we have: 
    
 \begin{equation} \label{FF}
   \begin{split}
   (2) = & \underbrace{ a_n \bP^{\ast}_{x_n}[   \mu_n \geq \eta n, S_n = y_n, \tau \geq \theta n]}_{(3)}\\ 
            & +  \underbrace{ a_n \bP^{\ast}_{x_n}[   \mu_n \geq \eta n, S_n = y_n, \tau < \theta n]}_{(4)}.
    \end{split}
 \end{equation}

  Making use of the Markov property, we have:
   \begin{equation} \begin{split}
    (3) & \leq a_n \sum_{j = 0}^{\gk a_n} \bP^{\ast}_{x_n}\left[ \max_{ i \leq \theta n} S_i \leq \gk a_n, S_{\theta n} = j, S_n = y_n\right] \\
         & \leq a_n \frac{V(y_n)}{V(x_n)} \sum_{j = 0}^{\gk a_n} \bP_{x_n}\left[ \max_{ i \leq \theta n} \hat{S}_i
 \leq \gk a_n, \hat{S}_{\theta n} = j\right] \bP_j\left[ \hat{S}_{(1- \theta)n} = y_n\right]\\
          & \leq a_n \frac{V(y_n)}{V(x_n)} \sum_{j = 0}^{\gk a_n} \bP_{x_n}\left[ \max_{ i \leq \theta n} S_i \leq \gk a_n, 
\tau_{ (-\infty,0)} > \theta n, S_{\theta n} = j\right] \bP\left[ S_{(1- \theta)n} = y_n - j\right],
   \end{split}
    \end{equation}
   where we recall that  $\tau_{ (-\infty,0)} = \inf \{ j \geq 1, S_j \in  (-\infty,0) \}$. 
    Using the local limit theorem and the fact that $j \in [0,\gk a_n]$, we get:
     
      \begin{equation} 
  \begin{split}
    & (3) \leq c \frac{V(y_n) \bP_{x_n}\left[ \tau_{ (-\infty,0)} > \theta n\right] }{V(x_n)} \\
 & \phantom{iiiiiiiiiiiiiiiiiiiiii} \times \bP_{x_n}\left[ \max_{ i \leq \theta n} S_i \leq \gk a_n \Big| \tau_{ (-\infty,0)} > \theta n\right] \frac{1}{\sqrt{1 - \theta}} 
 \phi \left( (\gep - \gk)(1 - \theta)^{-1/2} \right).
  \end{split}
             \end{equation} 
    Using the remark 4 in \cite{Shi}, we note that, as $n \to \infty$, 
     \begin{equation}
      \bP_{x_n}\left[ \max_{ i \leq \theta n} S_i \leq \gk a_n \Big| \tau_{ (-\infty,0)} > \theta n\right]
 \to m\left[ \sup_{[0,1]} \go_t \leq \frac{\gk}{\sqrt{\theta}} \right], 
       \end{equation}
   where $m(\cdot)$ denotes the measure of the brownian meander. 

    We prove that the equivalence
     \begin{equation}
        \bP_{x_n}[ \tau_{ (-\infty,0)} > \theta n] \sim V(x_n)   \bP[ T_1^- > \theta n] 
       \end{equation}
        holds uniformly on the sequences $x_n$ such that $x_n/a_n \to 0$ in Lemma \ref{LAA}, so that finally, using the convergence       
     \begin{equation}
      V(Ka_n)   \bP[ T_1^- > \theta n]  \to c \frac{K}{\sqrt{\theta}},
      \end{equation}
    which one can deduce from part \ref{FlT}, one gets:   
      \begin{equation} \begin{split}
    (3) & \leq c  V(Ka_n)   \bP[ T_1^- > \theta n] m \left[ \sup_{[0,1]} \go_t \leq \frac{\gk}{\sqrt{\theta}} \right]
 \frac{1}{\sqrt{1 - \theta}}  \phi \left( (\gep - \gk)(1 - \theta)^{-1/2} \right) \\
       & \leq c K m\left[ \sup_{[0,1]} \go_t \leq \frac{\gk}{\sqrt{\theta}}\right] \frac{1}{\sqrt{\theta(1 - \theta)}}
  \phi \left( (\gep - \gk)(1 - \theta)^{-1/2} \right)
    \end{split}
     \end{equation}
      and for $\theta > 0$ fixed, the quantity in the right hand side above vanishes as $\gk \searrow 0$. 

   We are left with the second term in equation \eqref{FF}. To get this, one notes that looking at the proof 
       of Lemma \ref{BBJ} in \cite{BJD},
    it is not difficult to see that, with $c,c' > 0$ fixed, the convergence in \eqref{CUU} 
          holds uniformly for $(u,v)$ in the compact set $[c,c'] \times [\gep,K]$.  Note in particular 
         the uniformity part in Lemma \ref{KL}, the fact that the convergence in the local limit
 theorem is uniform on the sets $[ca_n, c'a_n]$ and
       finally the fact that the derivative of the function $ (x,u) \mapsto \frac{x}{u^{3/2}} \phi( x/u^2)$ is uniformly bounded for
       $(x,u) \in [c,c'] \times (0,1)$ (to get the uniform 
    convergence of the Riemann's sums in the proof of Lemma \ref{BBJ} in \cite{BJD}).

           Making use once again of the Markov property, this implies that:   
          
           \begin{equation}
            \begin{split}
            (4) & \leq a_n \sum_{ j \leq \theta n} \sum_{k \geq \gk a_n} \bP^{\ast}_{x_n} [ \tau = j, S_j = k, \mu > \theta n] \bP^{\ast}_{k} [S_{n-j} = y_n]\\
              & \leq  \sum_{ j \leq \theta n} \sum_{k \geq \gk a_n} \bP^{\ast}_{x_n} [ \tau = j, S_j = k, \mu > \theta n]
 \frac{V(y_n)}{V(k)} a_n \bP_{k} \left[\hat{S}_{n-j} = y_n\right].
          \end{split}
            \end{equation}
            
             Note that one can restrict the range of summation of $k$ in the above expression over $[\gk a_n, K'a_n]$ where
    $K' > 0$ is large enough and independent of $n$. Thus, using Proposition \ref{BBJ}
  and the fact that $r(\cdot,\cdot)$ is continuous, one obtains:
             
              \begin{equation} \begin{split}
              (4) & \leq c \frac{  V(K)}{ V(\gk) \sqrt{1 - \theta} } \frac{K'}{\gep} \left[ \sup_{ u \in [\gk, K'], v \in [\gep, K]}
 r(u,v) \right] \sum_{ j \leq \theta n} \sum_{k \geq \gk a_n} \bP^{\ast}_{x_n} [ \tau = j, S_j = k, \mu_n > \theta n]\\
               & \leq  c \frac{  V(K)}{ V(\gk) \sqrt{1 - \theta} } \frac{K'}{\gep} \left[ \sup_{ u \in [\gk, K'], v \in [\gep, K]}
 r(u,v) \right] \bP^{\ast}_{x_n} [  \max_{ j \leq \mu_n} S_j \geq \gk a_n ]         
              \end{split}
               \end{equation}
                and as evidently the inclusion of events $\left\{  \max_{ j \leq \mu_n} S_j \geq \gk a_n  \right\} \subset
 \left\{ \max_{ j \leq \mu} S_j \geq \gk a_n \right\}$ holds,
  making use of Proposition \ref{BBBB}, the last term in the equation above vanishes as $n \to \infty$ since $x_n/a_n \to 0$.
   
  \subsection{ The case where $y_n/a_n$ vanishes at infinity} 
  This case relies heavily on the previous one. One has the equality:
 \begin{equation}
 \begin{split}
   & \bP_{x_n}\left[ \hat{S}_n = y_n\right] = \underbrace{\sum_{ z = \gep a_n}^{Ka_n} \bP_{x_n}\left[ \hat{S}_{n/2} = z\right]
 \bP_z \left[ \hat{S}_{n/2} = y_n\right]}_{(5)} \\
  & \phantom{iiiiiiiiiiiiiiiiiiiiiiiiiiiiiii}  +   \underbrace{\bP_{x_n}\left[ \hat{S}_n = y_n, S_{n/2} 
\leq \gep a_n, S_{n/2} \geq K a_n\right]}_{(6)}
 \end{split}
\end{equation}
 We first show that the term in (5) yields the desired estimate, and then that the term in (6) is 
negligible with respect to the first one. 

 For the term in $(5)$ , using time reversal and the case we just treated, one has of course:
  \begin{equation}
  \bP_z [ \hat{S}_{n/2} = y_n] \sim \frac{U(y_n) V(z)}{n/2} \bP[ S_{n/2} = z] 
\end{equation}
 so that, for $n \to \infty$, we have the equivalence:
  \begin{equation}
   \begin{split}
    \sum_{ z = \gep a_n}^{Ka_n} \bP_{x_n}\left[ \hat{S}_{n/2} = z\right] \bP_z \left[ \hat{S}_{n/2} = y_n\right]
   & \sim \sum_{ z = \gep a_n}^{Ka_n}  \frac{V(x_n) U(z)}{n/2}  \frac{U(y_n) V(z)}{n/2} \bP[ S_{n/2} = z]^2 \\  
   & \sim  \frac{V(x_n) U(y_n)}{n} \sum_{ z = \gep a_n}^{Ka_n} 8 \frac{z^2}{a_n^2} \frac{ \phi(z/a_{n/2})^2}{a_{n/2}^2} 
   \end{split}
\end{equation}
 where in the last equivalence we made use of the standard local limit theorem and of Lemma \ref{DDD}.
 Thus we are left with showing that 
  \begin{equation} \label{equi}
   \lim_{ \gep \searrow 0, K \nearrow \infty} \lim_{ n \to \infty} 8 \sqrt{2 \pi} a_n \sum_{ z = \gep a_n}^{Ka_n}
  \frac{z^2}{a_n^2} \frac{ \phi(z/a_{n/2})^2}{a_{n/2}^2} = 1. 
  \end{equation}
  We use Riemann's sum and the fact that $(a_n) \in \mathcal{R}_{1/2}$ to get:
   \begin{equation}
   \begin{split}
    8 \sum_{ z = \gep a_n}^{Ka_n}  \frac{z^2}{a_n^2} \frac{ \phi(z/a_{n/2})^2}{a_{n/2}^2} 
    & \sim 16 \sum_{ z = \gep a_n}^{Ka_n} \frac{z^2}{a_n^2} \frac{ \phi(\sqrt{2}z/a_n)^2}{a_n^2}\\
    & \sim \frac{16}{\sqrt{2 \pi}} \sum_{ z = \gep a_n}^{Ka_n} \frac{z^2}{a_n^2} \frac{ \phi(2z/a_n)}{a_n^2}\\
    & \sim \frac{16}{\sqrt{2 \pi} a_n}\int_{\gep}^K u^2 \phi(2u) du \\
    & \sim \frac{2}{\sqrt{2 \pi} a_n} \int_{\gep/2}^{K/2} u^2 \phi(u) du. 
   \end{split}
   \end{equation}
  
 and thus \eqref{equi} is valid. 

  We are left with showing that:
 \begin{equation} \begin{split}
               &  \limsup_{ K \nearrow \infty} \lim_{ n \to \infty} { n a_n \over V(x_n) U(y_n)} 
 \bP_{x_n}\left[ \hat{S}_n = y_n, \hat{S}_{n/2} \geq K a_n \right] = 0,  \\
               &   \limsup_{ \gep \searrow 0} \lim_{ n \to \infty} { n a_n \over V(x_n) U(y_n)}
  \bP_{x_n}\left[ \hat{S}_n = y_n, \hat{S}_{n/2} \leq \gep a_n \right] = 0.
                  \end{split}
    \end{equation}
    We define $\tilde{S}$ as being the time reversed version of $S$, that is the random walk whose transitions are given by
 \begin{equation}
   \bP[\tilde{S}_1 = y] := \bP[ S_1 = -y], \hspace{.6 cm} y \in \Z. 
 \end{equation}
    Note that 
   \begin{equation} \begin{split}
                   &  { n a_n \over V(x_n) U(y_n)} \bP_{x_n}\left[ \hat{S}_n = y_n, \hat{S}_{n/2} \geq K a_n\right]\\
                & = { n a_n \over V(x_n) U(y_n)} \sum_{z \geq Ka_n}  \bP_{x_n}\left[ \hat{S}_{n/2} = z\right]
 \bP_{y_n}\left[ \hat{\tilde{S}}_{n/2} = z\right].
                    \end{split}
    \end{equation}
  We recall that the following equivalences are shown in Lemma \ref{LAA} below:
 \begin{equation}
 \begin{split}
   & \bP_{x_n}[ \tau_{(- \infty, 0)} > n/2] \sim V(x_n) \bP[ T_1^- > n/2], \\
   &   \bP_{y_n}[ \tau_{(- \infty, 0)} > n/2] \sim V(y_n) \bP[ \tilde{T_1^-} > n/2]
 \end{split}
 \end{equation}
  and that they hold uniformly for $x_n, y_n$ which are $o(a_n)$. 

 Therefore, one deduces
 \begin{equation}
  \begin{split}
    &  { n a_n \over V(x_n) U(y_n)}  \bP_{x_n}[ \hat{S}_n = y_n, \hat{S}_{n/2} \geq K a_n ]  \\
    & \phantom{iiiii} \sim n  \bP[ T_1^- > n/2] \bP[ \tilde{T_1^-} > n/2] \\
   &  \phantom{iiiiiiiiiiiiii} \times  \sum_{z \geq Ka_n} a_n \bP_{x_n}[ S_{n/2} = z |
 \tau_{(- \infty, 0)} > n/2 ] \bP_{y_n}[ \tilde{S}_{n/2} = z | \tilde{\tau}_{(- \infty, 0)} > n/2 ].
  \end{split}
\end{equation}
   By the local limit theorem for the random walk conditioned to stay positive (see \cite[Theorem 2]{Car}):
 \begin{equation}
   \sup_{ z \in \Z} a_n \bP_{x_n}[ S_{n/2} = z | \tau_{(- \infty, 0)} > n/2 ] =: C < \infty.
 \end{equation}
  Recall that $T_1^+$ and $T_1^-$ are attracted to stable laws of index $1/2$, so that by standard Tauberian
 theorems (see \cite[ XIII 5.]{Fel2}):
 \begin{equation}
   \bP[ T_1^- > n] \sim \frac{1}{\sqrt{\pi}} \left( 1 - \bE\left[ e^{-\frac{1}{n}T_1^-}\right]\right), 
 \bP[ \tilde{T}_1^- > n] \sim \frac{1}{\sqrt{\pi}} \left( 1 - \bE\left[ e^{-\frac{1}{n} \tilde{T}_1^-}\right]\right).
 \end{equation}

   On the other hand, by the Wiener-Hopf factorization:
 \begin{equation}
   \begin{split}
     & 1 - \bE[ e^{- \gl T_1^-}] = \exp\left( - \sum_{n = 1}^{\infty} \frac{e^{- \gl n}}{n} \bP[ S_n < 0]\right) \\
    & 1 - \bE[ e^{- \gl T_1^+}] = \exp\left( - \sum_{n = 1}^{\infty} \frac{e^{- \gl n}}{n} \bP[ S_n \geq 0]\right)
   \end{split}
   \end{equation}
  hence, for $\gl \searrow 0$, 
 \begin{equation}
    \left(1 - \bE\left[ e^{- \gl T_1^-}\right]\right) \left( 1 - \bE\left[ e^{- \gl T_1^+}\right]\right) 
= \exp\left( - \sum_{n = 1}^{\infty} \frac{e^{- \gl n}}{n}\right) = 1 - e^{-\gl} \sim \gl
 \end{equation}
 therefore $\lim_{n \to \infty} n  \bP[ T_1^- > n]  \bP[ \tilde{T}_1^- > n] = \frac{1}{\pi}$. 
  Using finally the convergence towards the brownian meander, we get that
 \begin{equation}
   \begin{split}
        { n a_n \over V(x_n) U(y_n)}  \bP_{x_n}\left[ \hat{S}_n = y_n, \hat{S}_{n/2} \geq K a_n \right] & \leq
 \frac{C}{\pi} \bP_{y_n} \left[ \tilde{S}_{n/2} \geq K a_n \Big| \tilde{T}_1^- > n\right] \\
  & \leq c m\left[ \go_{1/2} > K\right]
  \end{split}
 \end{equation}
 and the last term vanishes as $K \to \infty$. Proceeding in the same way, it is easy to see that
  \begin{equation}
 \lim_{ n \to \infty} { n a_n \over V(x_n) U(y_n)}  \bP_{x_n}\left[ \hat{S}_n = y_n, \hat{S}_{n/2} \leq \gep a_n \right]
 \leq c m \left[\go_{1/2} \leq \gep \right],    
  \end{equation}
 and this last quantity also vanishes when $\gep \searrow 0$, and this concludes the proof of Lemma \ref{KL}.

 \section{Tightness of the measures $Q_n^{x_n,y_n}$} \label{SecTI}
 
  The proof of tightness is very similar to the one of \cite{BJD}.
  We first note that the process $S$ under $\bP^{\ast}_{x_n}[ \cdot | S_n = y_n]$ is still a Markov chain, so that according 
  to \cite[Theorem 8.4]{Bil} , tightness will follow if we can show that for each positive $\gep$ and
 $K \in (0,1)$, there exists $\gl > 0$ 
 and an integer $n_0$ such that
 \begin{equation} \label{TIT}
   \bP^{\ast}_{x_n}\left[ \max_{ i \leq Kn} S_i \geq \gl a_n \Big| S_n = y_n\right] \leq \frac{\gep}{\gl^2}
 \end{equation}
   for all $n \geq n_0 $.

   We proceed quite similarly as in the last part of the proof of Lemma \ref{KL}. We write:
   \begin{equation} \begin{split}
     & \star := \bP^{\ast}_{x_n}\left[ \max_{ i \leq Kn} S_i \geq \gl a_n \Big| S_n = y_n\right] \\
    & = \sum_{ j \geq 0} \bP_{x_n}\left[  \max_{ i \leq Kn} \hat{S}_i \geq \gl a_n, \hat{S}_{Kn} = j \Big| \hat{S}_n=y_n\right] \\
               & \sim \frac{ n a_n }{V(x_n)U(y_n)} \sum_{ j \geq 0} \bP_{x_n}
\left[  \max_{ i \leq Kn} \hat{S}_i \geq \gl a_n, \hat{S}_{Kn} = j\right] \bP_j\left[ \hat{S}_{n(1-K)} = y_n\right] \\                    
                            \end{split}
    \end{equation}
    Using the same considerations as in the last part of the proof of Lemma \ref{KL} (by simply replacing $n/2$ by $Kn$ or $(1-K)n$)
 , one gets that there exists a constant $C > 0$ such that: 
    \begin{equation}
     \star \leq   \frac{C}{ \sqrt{1-K}} \sum_{ j \geq 0} \bP_{x_n} \left[  \max_{ i \leq Kn} S_i \geq \gl a_n, S_{Kn} = j \Big|
 \tau_{ (- \infty,0)} > Kn \right] 
    \end{equation}
    so that using the weak convergence towards the brownian meander, we get: 
 \begin{equation}
   \star \leq \frac{C}{\sqrt{1-K}} m\left[ \sup_{t \in [0,1]} \go_t \geq \frac{\gl}{\sqrt{K}}\right],
 \end{equation}
  which for fixed $K$ vanishes exponentially fast when $\gl$ becomes large, and in particular \eqref{TIT} holds. This concludes the 
 proof of Theorem \ref{TA}, and thus we are done. 
   
    \section{Appendix} \label{SecV}

 The following is the main result of this appendix: 
  \begin{lemma} \label{LAA}
  Uniformly in $x_n$ such that $ x_n a_n^{-1} \to 0$ as $n \to \infty$, one has the following convergence:
  \begin{equation}\label{PF1}
    { \bP_{x_n}[\tau_{ (-\infty,0)} > n] \over \bP[T_1^- > n]} \sim V(x_n).
  \end{equation}
  \end{lemma}
    
  Note that it has been proved in \cite{BD} that
  \begin{equation}\label{BD}
    \liminf_{n \to \infty} { \bP_x[\tau_{ (-\infty,0)} > n] \over \bP[T_1^- > n]} \geq V(x)
  \end{equation}
 in full generality (that is for every oscillating random walk $S$ verifying $ \bP[S_1 > 0] \in (0,1)$).
   The convergence \eqref{PF1} has also been proved 
 in \cite{KOS} in the lattice case for fixed $x$. 
 \begin{proof}[Proof]
   For $x > 0$, we denote by $\tau_x = \inf\{ k \geq 1, S_k < -x \} $. One has the following identity:
  \begin{equation} \begin{split} 
     \bP_x[\tau_{ (-\infty,0)} > n] & = \bP[ \tau_x > n] \\
    & = \sum_{k = 0}^{ + \infty} \bP [ T^-_k \leq n < T^-_{k+1}, \tau_x > n] \\
     &= \sum_{k = 0}^{ + \infty} \bP [ T^-_k \leq n < T^-_{k+1}, H^-_k < x]\\
     &= \sum_{k = 0}^{ + \infty} \sum_{l=0}^n \bP [ T^-_k =l, H^-_k < x] \bP[ T^-_1 > n-l]
    \end{split}
 \end{equation}
 
 where in the last equality we made use of the Markov property. Thus:
  \begin{equation}\label{AS}
   \frac{\bP_x[\tau_{ (-\infty,0)} > n]}{V(x)} = \sum_{l=0}^n \frac{\bP[\hspace{.1 cm}  l \hspace{.1 cm} \text{is a descending ladder epoch}, -S_l< x]}{V(x)}  \bP[T^-_1 > n-l].  
  \end{equation}
   
 We recall a strong version of Iglehart's lemma (\cite[Lemma 5]{AliDon} ):
 \begin{lemma} \label{ALD}
   Let $c_n,d_n(z)$be two sequences  where $z$ belongs to a subset $\Delta$ of $\R$ . Define $e_n$ on $\Delta$ by:
 \begin{equation}
   e_n(z) := \sum_{j=0}^{n-1} d_j(z)c_{n-j}.
 \end{equation}
 Assume that there exist $c > 0$ such that the following condition holds uniformly on $z \in \Delta$:
   \begin{equation} \label{cdi1} 
     \sum_{j=1}^n d_j(z) \to d(z) < \infty \hspace{.2 cm} \text{and} \hspace{.2 cm} nd_n(z) \leq c \\
 \end{equation}
  Assume moreover that the sequence $c_n$ is regularly varying with index $ -\rho$ where $\rho \in (0,1)$. 
    Then the equivalence $e_n(z) \sim d(z) c_n $ holds uniformly on $ z \in \Delta$. 
 \end{lemma}
   We already pointed out that: 
 \begin{equation} \label{tT}
   \bP[T^-_1 > n ] \sim {b^{-1}(n) \over \sqrt{2 \pi}  n } \hspace{.2 cm} \text{as} \hspace{.2 cm} n \to \infty.
 \end{equation}
 Recalling that $b(\cdot) \in \mathcal{R}_{2}$, one has $b^{-1}(n)/n \in \mathcal{R}_{-1/2}$,
 which implies that the sequence  $\left(\bP[T^-_1 > n ]\right)_n$ 
 verifies the hypothesis of the sequence $c$ of Lemma \ref{ALD} with $\rho = 1/2$. 

   On the other hand, we write
  \begin{equation}
    1 = \sum_{ l \geq 0} \frac{\bP\left[\hspace{.1 cm}  l \hspace{.1 cm} \text{is a descending ladder epoch}, -S_l< x\right]}{V(x)} = \sum_{l \geq 0} \frac{ \sum_{ j \in [0, x]} v(l,j)}{V(x)}
  \end{equation}
  and thus we want to prove that the sequence $ d_l(x) =  \frac{\sum_{ j \in [0,x]} v(l,j)}{V(x)}$ satisfies the second conditions of 
 Lemma \ref{ALD} with $\gD_{(\gep_n)} = \left\{ (x_n) \in \Z^{\N}, \forall n, x_n \in [0, \gep_n a_n] \right\}$ where $\gep_n$ is a given positive 
  sequence which vanishes at infinity. 

  We first note that the uniform convergence of the series on $\gD_{(\gep_n)}$ has already been proved in 
the first part of the proof of Lemma \ref{KL}. 

  For the second point, we consider  a sequence $(x_n)_n \in \gD_{(\gep_n)}$. For $l > 0$ and making use of
 Lemma \ref{BJD1} (note in particular 
  the uniformity part of it) and of the local limit theorem:
  \begin{equation} \begin{split}
   \sum_{ j \in [0, x_n]} v(l,j) &  \leq x_n \sup_{ j \leq x_n} v(l,j)  \\
                                     & \leq  c \frac{x_n}{n a_n} V(x_n) \\
                                    & \leq c \frac{\gep_n}{n} V(x_n)
   \end{split}
  \end{equation}
  and as $\gep_n \to 0$,  both conditions of the first part
  of \eqref{cdi1} are fulfilled by the sequence
 $ \left( \frac{ \sum_{ j \in [0,x]} v(l,j)}{V(x)}\right)_l$.

  Thus we get that the following equivalence holds uniformly on $\gD_{(\gep_n)}$:
 \begin{equation}
       \frac{\bP_{x_n}\left[\tau_{(- \infty,0)} > n\right]}{V(x_n)} \sim \bP\left[T_1^- > n\right].
 \end{equation}
 This entails that the following equivalence holds uniformly on $x_n$ such that $x_n a_n^{-1} \to 0$ as $n \to \infty$:
 \begin{equation}
   \lim_{n \to \infty} { \bP_{x_n}\left[\tau_{(- \infty,0)} > n\right] \over \bP[T_1^- > n]} \sim V(x_n)
 \end{equation}
 which is equation \eqref{PF1}.
 \end{proof}     

 \textbf{Acknowledgement:} I am very grateful to Francesco Caravenna for his constant help and support during 
 this work. 
   
\bibliographystyle{alpha} 
\bibliography{Biblithese} 
\end{document}